\documentclass{amsart}

\usepackage{amsmath} 
\usepackage{verbatim}
\usepackage{textcomp}
\usepackage{amsfonts}
\usepackage{amssymb}
\usepackage{amsthm}
\usepackage{fancyhdr}
\usepackage{graphicx}
\usepackage[english]{babel}
\usepackage{wrapfig}
 \usepackage{caption} \usepackage{subcaption} 
\usepackage[margin=1in]{geometry}

\newcommand{\z}{\mathbb{Z}}

\newcommand{\bpf}{\begin{proof}}
\newcommand{\epf}{\end{proof}}
\newcommand{\ep}{\epsilon}

\newcommand{\iso}{\simeq}

\newtheorem{theorem}{Theorem}
\newtheorem{lemma}[theorem]{Lemma}

\newtheorem{corollary}[theorem]{Corollary}

\newenvironment{definition}[1][Definition]{\begin{trivlist}
\item[\hskip \labelsep {\bfseries #1}]}{\end{trivlist}}

\title{The graph structure of graph groups that are subgroups of Thompson's group $V$}
  
 \author{Nathan Corwin}
\address{Department of Mathematics - Hill Center 
Rutgers, The State University of New Jersey 
110 Frelinghuysen Rd. 
Piscataway, NJ 08854-8019}
\email{nacorwin@math.rutgers.edu}

\author{Kathryn Haymaker}
\address{Department of Mathematics and Statistics - St. Augustine Center, Villanova University
800 Lancaster Avenue 
Villanova, PA 19085
}
\email{kathryn.haymaker@villanova.edu }

\begin{document}

\begin{abstract}
We determine exactly which graph products, also known as Right Angled Artin Groups, embed into Richard Thompson's group $V$.
It was shown by Bleak and Salazar-Diaz that $\z^2 * \z$ was an obstruction. We show that this is the only obstruction.
This is shown by proving a graph theory result giving an alternate description of simple graphs without an appropriate induced subgraph.
\end{abstract}

\maketitle

	In this note we discuss which graph groups, also known as right angled Artin groups, exist as subgroups of Thompson's group $V$. 
	The group $V$ was first discovered in 1965 by the logician Richard Thompson while attempting to build algebras that encapsulate the properties of commutativity and associativity.
	Thompson's group $V$ quickly drew the interest of group theorists after Thompson showed that $V$ and a closely related group, Thompson's group $T$, were simple \cite{CFP96}.
	These were the first two known examples of finitely presented infinite simple groups.
	
	Despite interest in the group from various areas, the subgroup structure of $V$ is not well understood.
	Since $V$ contains many products of its subgroups, it was commonly believed that all graph groups embed in $V$. 
	This was disproved in 2013 by Bleak and Salazar-Diaz \cite{BOSD09} by showing $\z^2 * \z$ does not embed in $V$.
	Their motivation was related to computational complexity of the coword problem of groups and they only considered this particular graph group.
	The goal of this note is to classify exactly which  graph groups  embed in $V$.

	We interpret the question of which graph groups embed in $V$ as a graph theory question.
	In particular, we provide a useful description of finite simple graphs that do not have, as an induced subgraph, the graph associated to $\z^2 * \z$.
	We then show that graphs with this description correspond to graph groups that  embed in $V$.
	Therefore, all graph groups that do not contain $\z^2 * \z$  embed in $V$.
	Equivalently, the only graph groups that appear as subgroups of $V$ are direct products of free groups.
	
	The note is organized as follows. We give  background and definitions in Sections 1 and 2. 
	Section 1 discusses Thompson's group $V$ and graph groups.  Relevant  graph theory notions are explained in Section 2.
	We then formally state and prove our results in Section 3.
	
\section{Thompson's group $V$ and graph groups}

A graph group is a group that in a sense lies between free products and commutative free products. 
Given a finite simple graph $\Gamma$, the associated group $G(\Gamma) = G$ is a group whose generators are in bijective correspondence with the vertices of $\Gamma$.
The only relations are that some generators commute. 
In particular, two generators $x_i, x_j$ commute exactly when there is an edge between the corresponding vertices.

%

\begin{figure} 
	\centering{ 
		\begin{subfigure}[b]{0.4\textwidth} \includegraphics[scale=1]{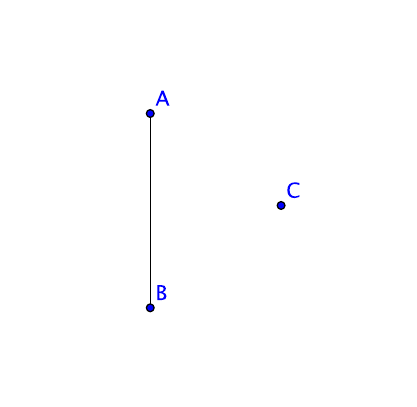} \caption{The graph associated with $\z^2*\z$} \label{ZfreeZsq} \end{subfigure}
		\quad
		\begin{subfigure}[b]{0.4\textwidth} \includegraphics[scale=1]{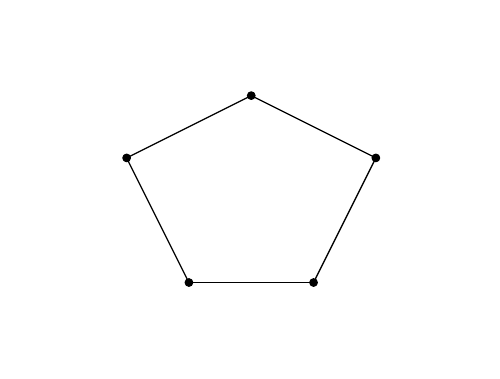} \caption{A Pentagon graph} \label{pentagon} \end{subfigure}		
	}
	\caption{Graph examples}	
\end{figure}	
 
A compete graph on $n$ vertices corresponds to $\z^n$ while the empty graph on $n$ vertices corresponds to the free group $F_n$.
As an intermediate example, consider the graph $\Gamma$ in Figure~\ref{ZfreeZsq}.
The corresponding graph group $G(\Gamma)$ will have a presentation $\langle a,b,c | ab=ba \rangle$, which factors as $\z^2 * \z$.
The graph of a pentagon shown in Figure~\ref{pentagon} represents a graph group that can not be described simply by using free products and direct products.

One common interpretation of Thompson's group $V$ is as a particular subgroup of the homeomorphism group of the Cantor set,
with elements of $V$ being the homeomorphisms that can be expressed as a finite list of prefix replacements.
It is difficult to give much more detail in a concise manner and the necessary facts about $V$ are listed below.
For more details, readers are referred to either 
\cite{CFP96} --- a standard introduction to Thompson's groups $F, T, $ and $ V$ --- 
or \cite{BOSD09} which has a very readable explanation of $V$ from a dynamical viewpoint.  
The following lemmas are well-known. 
\begin{lemma}
	\begin{enumerate}
		\item If $H_1$ and $H_2$ are subgroups of $V$, then there is a subgroup of $V$ isomorphic to $H_1 \times H_2$.
		\item There is a subgroup of $V$ isomorphic to the free group $F_2$.
	\end{enumerate}
\end{lemma}

As $\z$ is a subgroup of $F_2$, there is a subgroup of $V$ isomorphic to  $\z$. 
Thus, by repeatedly applying (1) we have that $\z^n$ is a subgroup of $V$ for all positive integers $n$.
Also, for any positive $n$, $F_n < F_2$, thus $F_n < V$.
Repeated applications of (1) show the following.

\begin{lemma}
	\label{DirectGood}
	If $G$ is a direct product of (finitely generated) free groups, then there is a subgroup of $V$ isomorphic to $G$.
\end{lemma}
Lastly, we state the theorem of Bleak and Salazar-D\'{i}az mentioned in the introduction. The following is Theorem 1.5 of \cite{BOSD09}.
\begin{theorem}
	\label{Collin}
	The group $\z^2 * \z$ does not embed in $V$.
\end{theorem}

\section{Relevant graph theory notions}

A finite simple graph $\Gamma$ is a set of vertices $V(\Gamma)$ and edges $E(\Gamma)$. $E(\Gamma)$ can be any subset of $V(\Gamma)\times V(\Gamma)$, where loops and multiedges are prohibited. That is, for any $v\in V(\Gamma)$,  $(v,v)\notin E(\Gamma)$, and an edge appears in $E(\Gamma)$ at most once. An induced subgraph of $\Gamma$ is a subgraph that is obtained by deleting vertices from $V(\Gamma)$.  
The eccentricity $\ep(v)$ of a vertex $v$ is the greatest geodesic distance between $v$ and any other vertex. It can be thought of as how far a vertex is from the vertex most distant from it in the graph.

We will be interested in two classes of graphs.
\begin{definition} Let $\mathcal{NB}$ be the class of all finite simple graphs $\Gamma$ that do not have the graph in Figure~\ref{ZfreeZsq} as an induced subgraph.
	In other words, $\mathcal{NB}$ is the class of finite simple graphs $\Gamma$ for which no  
	three distinct vertices $a,b,c \in \Gamma$ have the property that $(a,b) \in E(\Gamma)$ but both $(a,c),(b,c) \not\in E(\Gamma)$.
\end{definition}
	
\begin{definition}
	Denote by $\mathcal{GP}$ the class of all finite simple graphs $\Gamma$ for which there exists a partition 
	$P_0, P_1, \dots, P_n$,  of the vertices of $\Gamma$ with all, except perhaps $P_0$, nonempty such that
	\begin{itemize} 
	\item $P_0$ is the set of all vertices of $\Gamma$ with eccentricity one; 
	\item if $v_i, v_j \in P_k$ with $k > 0$, then $(v_i,v_j) \not\in E(\Gamma)$; 
	\item if $v_i \in P_i, v_j \in P_j$ with $i \ne j$ then $(v_i,v_j) \in E(\Gamma)$.
	\end{itemize} 
	We call an associated partition a \emph{commuting partition}.
\end{definition}
	In other words, a commuting partition separates the vertices of $\Gamma$ into the set $P_0$ of all vertices of eccentricity one and sets $P_1, \ldots, P_n$ that individually form empty subgraphs, but which have all possible edges between the subsets. 

\section{Results}

We are now prepared to prove the main theorem of this note.

\begin{theorem}
	\label{Same}
	The class $\mathcal{NB}$ is exactly the class $\mathcal{GP}$.
	
	\bpf
		First, suppose that $\Gamma$ is not in $\mathcal{NB}$. 
		Say $v_1, v_2, v_3 \in V(\Gamma)$ are distinct with $(v_1,v_2) \in E(\Gamma)$ but $(v_1, v_3),(v_2,v_3) \not\in E(\Gamma)$.
		Momentarily assume that $P_0, P_1, \dots P_n$ is a commuting partition for $\Gamma$.
		We have that $v_3 \in P_i$ with $i > 0$. 
		Since every vertex in $P_i$ is adjacent to every vertex  in $P_j$ where $j\not=i$, and $v_3$ is not adjacent to $v_1$, we must have that $v_1 \in P_i$.
		Similarly, since $v_3$ and $v_2$ are not adjacent, $v_2 \in P_i$. 
		Thus, $v_1, v_2 \in P_i$ implies $(v_1,v_2) \not\in E(\Gamma)$ but we are assuming $(v_1,v_2) \in E(\Gamma)$.
		This contradiction shows that there can not be a commuting partition for $\Gamma$ thus $\Gamma$ is not in $\mathcal{GP}$.
		Equivalently, $\mathcal{GP} \subset \mathcal{NB}$.
		
		Now, assume that $\Gamma \in \mathcal{NB}$. We will give an algorithm to partition the vertices, then we will show the resulting partition is a commuting partition.
		
		For the initial step, set $P_0 = \{ v \in V(\Gamma) | \ep(v) = 1\}$.
		For the recursive step, assume that $P_0, P_1, \dots, P_k$ is a partial partition of $\Gamma$.
		Let $S_k = \cup_{i=0}^k P_i$ and $R_k = V(\Gamma) \setminus S_k$.
		If $S_k = V(\Gamma)$, then the already built $P_i$'s are a partition of $\Gamma$ and we stop.
		Otherwise, choose $w_k \in R_k$.
		Set $P_{k+1}= \{ v \in R_k | (w_k,v) \not\in E(\Gamma) \}$. 
		As $\Gamma$ is simple, $w_k \in P_{k+1}$.
		Therefore, the size of $R_{k+1}$ is always at least one smaller than $R_k$ and, as $V(\Gamma)$ is finite, this process must terminate.
		
		We now show the resulting partition is a commuting partition. By construction, $P_0$ is exactly as required and $P_i$ is nonempty if $i >0$.
		
		For condition (ii), consider $v_i,v_j \in P_{k+1}$. 
		If either is $w_k$ then by construction they are not adjacent.
		Otherwise neither are adjacent to $w_k$ and thus can't be adjacent to each other as $\Gamma$ is in $\mathcal{NB}$ and this would form an inadmissible induced subgraph.
		
		For condition (iii), if  $v_i \in P_i, v_j \in P_j$ with $i < j$, we again consider cases.
		If $v_i = w_{i-1}$, then we conclude that $(v_i,v_j)\in E(\Gamma)$ as this is the only reason $v_j$ would have been excluded from $P_i$.
		Otherwise, we have $(v_i, w_{i-1}) \not\in E(\Gamma)$ and $(v_j,w_{i-1}) \in E(\Gamma)$.
		As $\Gamma$ is in $\mathcal{NB}$, we must have $(v_i,v_j) \in E(\Gamma)$ or else $v_i,v_j,w_{i-1}$ forms an  inadmissible induced subgraph.
		
		We conclude that the partition is a commuting partition hence $\Gamma$ is in $\mathcal{GP}$.
				
	\epf
\end{theorem}

\begin{lemma}
	\label{KGgood}
	Let $\Gamma$ be in $\mathcal{GP}$. Then, the corresponding group $G(\Gamma)$ is a subgroup of $V$.
	
	\bpf
		Let $\mathcal{P} $ be a commuting partition for $\Gamma$. 
		For each $i$, consider the induced subgraph $\Gamma_i$ formed by looking at only the vertices in $P_i$.
		The induced subgraph $\Gamma_0$ corresponds to a group $G_0 \iso \z^{|P_0|}$.
		For each positive $i$, the induced subgraph $\Gamma_i$ corresponds to a group $G_i \iso F_{|P_i|}$.
		As vertices from different elements of the partition have an edge between them, the corresponding generators of $G(\Gamma)$ commute.
		Thus, $G(\Gamma) \iso G_0 \times G_1 \times \cdots  \times G_k$.
		This is a subgroup of $V$ by Lemma~\ref{DirectGood}.
	\epf
\end{lemma}

\begin{corollary}
	Each of the following is true:
	\begin{itemize}
		\item A graph group embeds into $V$ if and only if its associated graph is in $\mathcal{GP}$.
		\item A graph group embeds into $V$ if and only if it is the direct product of free groups.
		\item $\z^2 * \z$ is the only obstruction to a graph group embedding into $V$.
	\end{itemize}

	\bpf
		Lemma~\ref{KGgood} proved that if $\Gamma$ is in $\mathcal{GP}$, then $G(\Gamma)$ embeds into $V$. 
		Theorem~\ref{Collin} implies that if $\Gamma$ is not in $\mathcal{NB}$ then $G(\Gamma)$ does not embed into $V$ and
		Theorem~\ref{Same} shows that $\mathcal{NB}$ is equivalent to $\mathcal{GP}$.
		Therefore, if $\Gamma$ is not in $\mathcal{GP}$ then $G(\Gamma)$ does not embed into $V$.
	\epf
\end{corollary}

\bibliographystyle{plain}
\bibliography{KandN}

\end{document}